\documentclass[10pt,twoside,reqno]{amsart}
\usepackage{amssymb}
\calclayout

\numberwithin{equation}{section}
\makeatletter

\renewcommand{\@secnumfont}{\bfseries}

\renewcommand{\section}{\@startsection{section}{1}%
  {0mm}{.7\linespacing\@plus\linespacing}{.5\linespacing}
  {\normalfont\bfseries\centering}}

\newcommand{\bibsection}{\@startsection{section}{1}%
  {0mm}{.7\linespacing\@plus\linespacing}{.5\linespacing}
  {\normalfont\scshape\centering}}

\renewcommand{\@biblabel}[1]{#1.}

\newtheorem{thm}{\bf Theorem}[section]

\begin{document}

\vspace{1.3cm}

\title
        {Some identities of $q$-Bernoulli  numbers associated $p$-adic convolutions}

\author{J.J. Seo, T.Kim, S.H.Lee}

\thanks{\scriptsize }

\address{1\\Department of Applied Mathematics\\
            Pukyong National University\\
            Busan 608-737, Republic of Korea}
\email{seo2011@pknu.ac.kr}

\address{2\\Department of  Mathematics\\
            Kwangwoon National University\\
           Seoul 139-701, Republic of Korea}
\email{tkkim@kw.ac.kr}
\address{3\\Division of  General Education\\
            Kwangwoon National University\\
           Seoul 139-701, Republic of Korea}
\email{leesh58@kw.ac.kr}

\subjclass{  bbbbb}

\maketitle

\begin{abstract} In this paper, we give some interesting and new identities of $q$-Bernoulli numbers which are derived from convolutions on the ring of $p$-adic integers.

\end{abstract}

\pagestyle{myheadings}
\markboth{\centerline{\scriptsize J.J. Seo, T.Kim, S.H.Lee}}
          {\centerline{\scriptsize Some identities of $q$-Bernoulli numbers associated $p$-adic convolutions}}
\bigskip
\bigskip
\medskip
\section{\bf Introduction}
\bigskip
\medskip

Throughout this paper $ \mathbb{Z}_p, \mathbb{Q}_p$ and $\mathbb{C}_p$ will denote the ring of $p$-adic integers, the field of $p$-adic numbers and the completion of algebraic closure of $\mathbb{Q}_p$. The $p$-adic norm is defined by $|p|_p=p^{-v_p(p)} =p^{-1}$. Now, we set
\begin{equation*}
U_p=\left\{\alpha\in \mathbb{C}_p | \ |\alpha-1|_p <1 \right\}
\end{equation*}
and
\begin{equation*}
T_p=\left\{w\in \mathbb{C}_p | w^{p^n} =1 \ \ \textnormal{for some}\ \ n\geq 0 \right\},
\end{equation*}
so that $T_p$ is the union of cyclic (multiplicative) group $\mathbb{C}_{p^n}$ of order $p^n$ $(n\geq 0)$ and $T_p\subset U_p$ (see [12]).
Let $UD(\mathbb{Z}_p)$ and $C(\mathbb{Z}_p)$ be the space of uniformly differentiable and continuous function on $\mathbb{Z}_p$.
For $f\in UD(\mathbb{Z}_p)$, the $p$-adic invariant integral on $\mathbb{Z}_p$ is defined as follows :
\begin{equation}
I_0(f)=\int_{\mathbb{Z}_p} f(x) d\mu_0 (x) = \lim_{N\rightarrow \infty} {1\over {p^N}} \sum_{x=0}^{p^N-1} f(x), \  (\textnormal{see}\ [11]).
\end{equation}

For $w \in T_p$, the $p$-adic Fourier transform of $f$ is given by
\begin{equation}
{\hat f}_w = I_0 (f\phi_w )=\lim_{N\rightarrow \infty} {1\over {p^N}} \sum_{x=0}^{p^N-1} f(x)\phi_w (x),  \  (\textnormal{see}\ [12]),
\end{equation}
where $\phi_w(x)=w^x.$

Let $f, g\in UD(\mathbb{Z}_p).$ Then \textit{C. F. Woodcock} defined the convolution of $f$ and $g$ as follows :
\begin{equation}
f*g=\sum_w {\widehat f}_w {\widehat g}_w \phi_{w^{-1}},   \ \ (\textnormal{see}\ [12]).
\end{equation}

From (1.3), we note that
\begin{equation}
f*g \in UD(\mathbb{Z}_p),\ \ (\widehat{f*g})_w ={\widehat f}_w {\widehat g}_w, \ (\forall w\in T_p),
\end{equation}
and $(UD(\mathbb{Z}_p), +, *, V)$ is Banach algebra, where $V(f)=\min\left\{\nu(f), R(f) \right\}, \nu(f)=\inf_{x\in\mathbb{Z}_p} \nu_p\left(f(x) \right)$ and $R(f)=\inf \nu_p \left({{f(x)-f(y)}\over {x-y}} \right), x\not=y \in \mathbb{Z}_p$,  (see [12]).

Let $int$ $\mathbb{Z}_p=\left\{f\in UD(\mathbb{Z}_p) | f^{\prime}=0 \right\}$. Then $int$ $\mathbb{Z}_p$ is $*$-ideal of $UD(\mathbb{Z}_p)$. Differentiation induced a linear isometry
\begin{equation}
UD(\mathbb{Z}_p)/ {int } \\ \mathbb{Z}_p  \longrightarrow C(\mathbb{Z}_p) \ \ \textnormal{by}\ \ (f*g)^{\prime}=-f^{\prime} \circledast g^{\prime},
\end{equation}
where $f, g \in UD(\mathbb{Z}_p) $ (see [12]). By (1.5), we get
\begin{equation}
(f\circledast g)(n)=\sum_{i=0}^n f(i)g(n-i),
\end{equation}
where $f, g \in C(\mathbb{Z}_p).$

For $f, g \in UD(\mathbb{Z}_p)$ and $z \in \mathbb{Z}_p$, it is known that
\begin{equation}
f*g(z)=I_0^{(x)} (f(x)g(z-x))-f \circledast g^{\prime}(z),
\end{equation}
where $I_0^{(x)}(f) = \lim_{N\rightarrow \infty}  \frac{1}{p^N}\sum_{x=0}^{p^N -1} f(x)$.

When one talks of $q$-extensions, $q$ is variously considered as an indeterminate, a complex number $q\in \mathbb{C}$, or a $p$-adic number $q\in \mathbb{C}_p$.
If $q\in \mathbb{C}_p$, then we normally assume $|1-q|_p < p^{- 1\over {p-1}}$ so that  $q^x=\exp(\log q) $ for $|x|_p \leq 1$. We use the notation $\left[x\right]_q=\left[x:q\right]={{1-q^x}\over{1-q}}$. Thus, $\lim_{q\rightarrow 1}\left[x\right]_q =x$.

As is well known, the usual \emph{$Bernoulli \ numbers$} are defined by
\begin{equation}
B_0=1, \ \ (B+1)^n -B_n =
\left\{
\begin{array}{cc}
1 , & \text{if }n=1 \\
0, & \text{if\ }n> 1\text{,}%
\end{array}%
\right.
\end{equation}

with the usual convention about replacing $\beta^i$ by $\beta_i$ (see [1-10]).
In [1.3], $Carlitz$ defined the $q$-extension of $Bernoulli\ numbers$ as follows:

\begin{equation}
\beta_{0,q} =1, \ \ q (q\beta_q+1)^n -\beta_{n,q}=
\left\{
\begin{array}{cc}
1 , & \text{if }n=1 \\
0, & \text{if\ }n> 1\text{,}%
\end{array}%
\right.
\end{equation}

with the usual convention about replacing $\beta_q^i$ by $\beta_{i,q}$.
By (1.9), we easily see that
\begin{equation*}
\beta_{0,q}=1,\ \beta_{1,q}={-1\over {[2]_q}},\ \beta_{2,q}={q\over {[2]_q [3]_q}},\ \beta_{3,q}={{(1-q)}\over {[3]_q [4]_q}}, \cdots.
\end{equation*}

In this paper, we consider the modified $q\emph{-Bernoulli\ numbers},$ which are slightly different $\emph{Carlitz's}$ $q$\emph{-Bernoulli numbers}, as follows :
\begin{equation}
\widetilde{\beta}_{0,q} =1, \ \left(q\widetilde{\beta}_q+1 \right)^n -\widetilde{\beta}_{n,q}=
\left\{
\begin{array}{cc}
\frac{logq}{q-1} , & \text{if }n=1 \\
0, & \text{if\ }n> 1\text{,}%
\end{array}%
\right.
\end{equation}
with the usual convention about replacing $\widetilde{\beta}_q^i$ by $\widetilde{\beta}_{i,q}.$

The purpose of our paper is to give some interesting and new identities of the modified $q-Bernoulli\ numbers$ $\widetilde{\beta}_{n,q}$ which are derived from convolutions on the ring of $p$-adic integers.
\bigskip
\bigskip

\section{\bf Some identities of $q$-Bernoulli  numbers}
\bigskip
\medskip
Let us consider the following $q$-extension of $Bernoulli$  $numbers$ :
\begin{equation}
\widetilde{\beta}_{0,q}=1, \ \ \left(q\widetilde{\beta}_q+1\right)^n-\widetilde{\beta}_{n,q} =
\left\{
\begin{array}{cc}
\frac{logq}{q-1} , & \text{if }n=1 \\
0, & \text{if\ }n> 1\text{,}%
\end{array}%
\right.
\end{equation}

Then, by (1.1), we easily see that

\begin{equation}\begin{split}
\widetilde{\beta}_{n,q}=&\int_{\mathbb{Z}_p} [x]_q^n d\mu_0(x)\\
                       =&{{\log q}\over{(1-q)^{n+1}}}\sum_{l=0}^n \binom{n}{l} (-1)^{l-1} {l \over{[l]_q}}, \ \ n\in \mathbb{Z}_+ = \mathbb{N}\cup \left\{0\right\}.
\end{split}\end{equation}

In the equation (1.6) and (1.7), if we take $p$-adic integral on $\mathbb{Z}_p$ with respect to variable $\mathbb{Z}$, then we have
\begin{equation}
I_0^{(z)}(f*g)=I_0^{(z)}\left( I_0^{(x)}(f(x)g(z-x))\right)-I_0^{(z)} \left(f\circledast g^{\prime} (z) \right).
\end{equation}

By (1.4) and (2.3), we get

\begin{equation}
I_0^{(z)}\left(f\circledast g^{\prime}(z)\right)=I_0^{(z)}\left( I_0^{(x)}(f(x)g(z-x))\right)-I_0^{(z)} (f) I_0^{(z)}(g).
\end{equation}
Let us take $f(x)=[x]_{q^{-1}}^m,$ $g[x]=[x]_q^n$.
Then
\begin{equation*}
q^{\prime}(x)=n{{\log q}\over{q-1}} [x]_q^{n-1}, \ \ (m, n \in \mathbb{N}).
\end{equation*}

Now, we set
\begin{equation}
A_{m,n}^q =I_0^{(z)}\left([z]_{q^{-1}}^m \circledast[z]_q^{n-1} \right), \ \ m, n \in \mathbb{N}.
\end{equation}

From (2.4) and (2.5), we can derive

\begin{equation}
n{{\log q}\over{q-1}}A_{m,n-1}^q = \int_{\mathbb{Z}_p}\int_{\mathbb{Z}_p} [x]_{q^{-1}}^m [z-x]_q^{n} d\mu_0(x) d\mu_0(z)-\widetilde{\beta}_{m,q^{-1}} \widetilde{\beta}_{n,q}.
\end{equation}
Note that
\begin{equation}
[z-x]_{q^{-1}}^n=\left([z]_q-q^{-1}q^z [x]_{q^{-1}} \right)^n =\sum_{l=0}^n  \binom{n}{l}   [z]_q^{n-l}(-1)^lq^{-l}q^{lz} [x]_{q^{-1}}^l.
\end{equation}

By (2.6) and (2.7), we get

\begin{equation}\begin{split}
n{{\log q}\over{q-1}}A_{m,n}^q =&\sum_{l=1}^n\binom{n}{l} (-1)^l q^{-l}   \int_{\mathbb{Z}_p} [x]_{q-1}^{m+l} d\mu_0(x)\int_{\mathbb{Z}_p} [z]_q^{n-l} q^{lz} d\mu_0 (z) \\
                               =&\sum_{l=1}^n \binom{n}{l}(-1)^k q^{-l} \widetilde{\beta}_{m+l,q^{-1}} \sum_{k=0}^l \binom{l}{k}(q-1)^k \int_{\mathbb{Z}_p} [z]_q^{n+k-l}d\mu_0 (z)\\
                               =&\sum_{l=1}^n\sum_{k=0}^{l} \binom{n}{l}\binom{l}{k} (-1)^l q^{-l} \widetilde{\beta}_{m+l,q^{-1}} (q-1)^k \widetilde{\beta}_{m+k-l,q}.
\end{split}\end{equation}

Thus, from (2.8), we have

\begin{equation}
A_{m,n}^q={{q-1}\over{\log q}}{1\over n}\sum_{l=1}^n \sum_{k=0}^l \binom{n}{l} \binom{l}{k}(-1)^lq^{-l}(q-1)^k \widetilde{\beta}_{m+l,q^{-1}}\widetilde{\beta}_{n+k-l,q}.
\end{equation}

By (2.5), we easily see that

\begin{equation}\begin{split}
A_{m,n}^q=&I_0^{(z)} \left([z]_{q^{-1}}^m\circledast [z]_q^{n-1} \right)\\
         =&I_0^{(z)} \left([z]_q^{n-1}\circledast [z]_{q^{-1}}^m \right)\\
         =&A_{n-1,m+1}^q.
\end{split}\end{equation}

From (2.5), we have

\begin{equation}\begin{split}
\nu_p(A_{m,n}^q)=&\nu_p \left(I_0\left([z]_{q^{-1}}^m\circledast [z]_q^{n-1} \right) \right)\\
                \geq &\nu \left([z]_{q^{-1}}^m\circledast [z]_q^{n-1} \right)-1 \\
                \geq &\nu ([z]_{q^{-1}}^m)+\nu([z]_q^{n-1})-1-1\\
                \geq &-2.
\end{split}\end{equation}

Therefore, by (2.9), (2.10) and (2.11), we obtain the following theorem.

\begin{thm}\label{thm 2.1} For $m\in \mathbb{Z}_+$ and $n \in \mathbb{N}$, we have
\begin{equation*}
A_{m,n}^q={{q-1}\over{\log q}} {1\over n} \sum_{l=1}^n\sum_{k=0}^l  \binom{n}{l}\binom{l}{k} (-1)^lq^{-l}(q-1)^k \widetilde{\beta}_{m+l,q^{-1}}\widetilde{\beta}_{n+k-l,q} .
\end{equation*}
Furthermore,
\begin{equation*}
A_{m,n}^q=A_{n-1,m+1}^q, \ \ \nu_p(A_{m,n}^q)\geq -2.
\end{equation*}
\end{thm}

In particular, if we take $m=0$, then by (2.10), we get

\begin{equation}
A_{0,n}^q=A_{n-1,1}^q \ \ \textnormal{($q$-analogue of \textit{Euler identity})}.
\end{equation}

From (2.2), we note that

\begin{equation}\begin{split}
\widetilde{\beta}_{n,q}=&{{\log q}\over{(1-q)^{n+1}}}\sum_{l=0}^n\binom{n}{l} (-1)^{l-1} l {{1-q}\over{1-q^l}}\\
                       =&n{{\log q}\over{(1-q)^n}}\sum_{l=1}^n\binom{n-1}{l-1} (-1)^{l-1} \sum_{m=0}^\infty q^{lm}\\
                       =&n{{\log q}\over{1-q}}{1\over{(1-q)^{n-1}}} \sum_{l=0}^{n-1} \binom{n-1}{l}(-1)^{l} \sum_{m=0}^\infty q^{(l+1)m}\\
                       =&n{{\log q}\over{1-q}} \sum_{m=0}^{\infty} q^m [m]_q^{n-1},
\end{split}\end{equation}

where $n \in \mathbb{N}$. Thus, by (2.13), we obtain the following theorem.

\begin{thm}\label{thm 2.2} For $n\in \mathbb{N}$, we have
\begin{equation*}
- {{\widetilde{\beta}_{n,q}}\over{n}} ={{\log q}\over{q-1}}\sum_{m=1}^\infty q^m [m]_q^{n-1}.
\end{equation*}
\end{thm}

Let $F_q(t)$ be the generating function for $\widetilde{\beta}_{n,q}$ with $F_q(t)=\sum_{k=0}^\infty {1\over{k!}}\widetilde{\beta}_{k,q}t^k.$
Then, by (2.13), we get

\begin{equation}\begin{split}
F_q(t)=&\sum_{k=0}^\infty {\widetilde{\beta}_{k,q}\over{k!}} t^k \\
      =&{{\log q}\over{1-q}}\sum_{m=0}^{\infty} q^m \left\{1+\sum_{k=1}^\infty {k\over{k!}} [m]_q^{k-1}t^k \right\}\\
      =&{{\log q}\over {(1-q)^2}}+t{{\log q}\over{1-q}}\sum_{m=0}^\infty q^m e^{[m]_q^t}.
\end{split}\end{equation}
Therefore, by (2.14), we obtain the following theorem.

\begin{thm}\label{thm 2.3} Let $F_q(t)=\sum_{k=0}^\infty \widetilde{\beta}_{k,q} {t^k \over {k!}}$. Then we have
\begin{equation*}
F_q(t)={{\log q}\over{(1-q)^2}} +t{{\log q}\over{1-q}}\sum_{m=0}^\infty q^m e^{[m]_q^t}.
\end{equation*}
\end{thm}

\end{document}